\documentclass[12pt,letterpaper,reqno]{amsart}

\usepackage{times}
\usepackage[T1]{fontenc}
\usepackage{mathrsfs}
\usepackage{latexsym}
\usepackage[dvips]{graphics}
\usepackage{epsfig}
\usepackage{amsmath,amsfonts,amsthm,amssymb,amscd}
\input amssym.def
\input amssym.tex

\newcommand\be{\begin{equation}}
\newcommand\ee{\end{equation}}
\newcommand\bea{\begin{eqnarray}}
\newcommand\eea{\end{eqnarray}}
\newcommand\bi{\begin{itemize}}
\newcommand\ei{\end{itemize}}
\newcommand\ben{\begin{enumerate}}
\newcommand\een{\end{enumerate}}
\newcommand\bc{\begin{center}}
\newcommand\ec{\end{center}}
\newcommand\ba{\begin{array}}
\newcommand\ea{\end{array}}

\newtheorem{thm}{Theorem}[section]

\newtheorem{cor}[thm]{Corollary}
\newtheorem{lem}[thm]{Lemma}

\theoremstyle{definition}

\begin{document}

\title{Minimal odd order automorphism groups}

\author{Peter Hegarty and Desmond MacHale}
\email{hegarty@math.chalmers.se, d.machale@ucc.ie} \address{Mathematical Sciences,
Chalmers University
Of Technology and University of Gothenburg, 41296 G\"oteborg, Sweden}
\address{Department of Mathematics, University College Cork, Cork, Ireland}

\date{\today}

\begin{abstract} We show that $3^7$ is the smallest order of a non-trivial 
odd order group
which occurs as the full
automorphism group of a finite group.
\end{abstract}


\maketitle

\setcounter{equation}{0}

\setcounter{equation}{0}

\section{Introduction and notation}

The map which takes every element of a group 
to its inverse is an automorphism
if and only if the group is abelian. Following \cite{HL}, we say that a
group is N.I. (no inversions occur) if no automorphism of the group sends
any non-identity element to its inverse{\footnote{In \cite{MS}, the term 
I.F.P. (inverse-point-free) is used instead.}}. Thus an N.I. group is 
non-abelian. It is easy to see that a finite group $G$ is N.I. if and
only if Aut$(G)$ has odd order{\footnote{We do not know if there exists an 
infinite group which has no automorphism of order $2$ but is not N.I.}}, and 
that the order of such a group must be 
odd. N.I. groups have attracted some attention, because of the
surprising difficulty in finding explicit
examples of them. Though there are a number of results indicating that
such groups are ubiquitous (for example, see \cite{HL}, 
\cite{HR} and \cite{M}), no brute force search among $\lq$small' groups
will yield any examples. More precisely,  
MacHale and Sheehy \cite{MS} proved that there is no finite N.I.
group of order less than $3^6$. There
exist groups of order $3^6$ with automorphism group of order $3^7$ (see
the Small Groups Library \cite{BEB}). In this
paper we prove the following complement to MacHale and Sheehy's result :   

\begin{thm}
If $G$ is a non-trivial finite group such that 
$|{\hbox{Aut}}(G)|$ is odd, then \\ $|{\hbox{Aut}}(G)| \geq 3^7$.
\end{thm}

The proof will be presented in the next two sections. In Section 4 we
present some remaining open problems concerning small N.I. groups. 
\\
\\
All groups considered in the remainder of this paper are finite. The 
following notational conventions will be adhered to throughout.
\\
\\
We denote by $\Omega(n)$ the number of prime factors of the
positive integer $n$, counting repetitions, and by $\omega(n)$ the number of distinct prime factors of $n$. If $G$ is a group then $\pi(G)$ denotes
the set of distinct prime
factors of $|G|$. If $p \in \pi(G)$ then $G_p$ denotes a Sylow $p$-subgroup of $G$.
\par
The direct product of $n$ copies of a group $G$ with itself will be denoted $G^n$.
\par
The center of a group $G$ is denoted by 
$Z(G)$. We will often just write $Z$ when it is clear
to which group we are referring.
\par
Inn$(G)$ and Cent$(G)$ denote, respectively, the groups of inner and central
automorphisms
of $G$. For $g \in G$, the inner automorphism $x \mapsto g^{-1}xg$ will be denoted
$I_g$. The commutator $x^{-1}g^{-1}xg$ will be denoted $[x,g]$.
\par
Let $A$ be a subgroup of $G$. Let $\phi
\in {\hbox{Aut}}(A)$ and $\phi^{*} \in {\hbox{Aut}}(G)$. We say
that $\phi^{*}$ is an {\em extension} of $\phi$ (or, equivalently, that $\phi$ {\em
extends} to $\phi^{*}$), if
$(A)\phi^{*} = A$ and $\phi^{*} |_{A} = \phi$.
\par
If $A$ is a normal subgroup of $G$, we say that $G$ {\em splits} over $A$ if there
is a subgroup $B \subseteq G$ such that $A \cap B = \{1\}$ and $G = AB$. Then $G$ is
said to be a {\em split extension} of $A$ by $B$ and we write $G = A \rtimes B$. The
subgroup $B$ is said to be a {\em complement} of $A$. It acts on $A$ by conjugation
and the induced subgroup of Aut$(A)$ is denoted by 
$A^{B}$. This last piece of notation
will be used even in the case of a non-split extension.
\par If $q$ is a prime power, then $\mathbb{F}_q$ denotes the finite field with $q$
elements and
GL$(n,q)$ the group of invertible linear transformations on the
$n$-dimensional vector space over
$\mathbb{F}_q$. If $q$ is a prime, this group can be identified with Aut$(G)$ where
$G \cong C_{p}^{n}$, the
elementary abelian $p$-group of rank $n$.
\par

\setcounter{equation}{0}

\section{Toolbox}

All of the results in this section, except for Theorem 2.13, are taken from the existing literature. The first seven lemmas are, moreover, standard material. 
For these the reader seeking 
details beyond what we provide 
is referred to \cite{Sco} in the absence of any alternative citations.

\begin{lem}
Cent$(G) \; \cap$ Inn$(G) = Z[{\hbox{Inn}}(G)] \cong Z[G/Z(G)]$.
\end{lem}


\begin{lem}
Let $G$ be a group without any abelian direct factor. Then
\be
|{\hbox{Cent}}(G)| = |{\hbox{Hom}}(G/G^{\prime},Z(G))|.
\ee
\end{lem}

\begin{proof}
The natural one-to-one correspondence is given by $\phi \rightarrow
(G^{\prime}g \mapsto g^{-1} g\phi)$. For a rigorous proof of the lemma, see
\cite{Sa}.
\end{proof}

\begin{lem}
Let $G$ be a group with center $Z$ and suppose $G/Z$ is a $p$-group for some prime
$p$.
\\
(i) If $G/Z$ is abelian, say
\be
G/Z \cong \prod_{i=1}^{r} C_{p^{\alpha_i}},
\ee
with $\alpha_1 \geq \alpha_2 \geq \cdots \geq \alpha_r$, then $r \geq 2$ and
$\alpha_1 = \alpha_2$.
\\
(ii) If $G/Z$ is not elementary abelian then it possesses a non-trivial central
characteristic subgroup $K$
of index at least $p^2$.
\end{lem}

\begin{proof}
(i) Clearly $r \geq 2$.
Let $Za$ be a generator of a highest order cyclic factor. Then for any $Zb$ in
any complementary factor we
have $b^{p^{\alpha_2}} \in Z$. Since $G/Z$ is abelian we have $1 =
[a,b^{p^{\alpha_2}}]= [a,b]^{p^{\alpha_2}} = [a^{p^{\alpha_2}},b]$.
Thus $a^{p^{\alpha_2}} \in Z$, so $\alpha_1 = \alpha_2$.
\\
(ii) If $G/Z$ is not abelian, take $K$ to be its center. Otherwise take $K$ to be
the subgroup generated by
the elements of order $p$. In both cases $K$ is non-trivial of index at least $p^2$,
by part (i).
\end{proof}

\begin{lem}
Let $G$ be a group, $N$ a normal subgroup of $G$. Let $\phi$ be an automorphism of
$G$ whose order is
relatively prime to that of $G$ and which acts trivially on both $N$ and $G/N$. Then
$\phi$ is the trivial
automorphism.
\end{lem}

\begin{proof}
Let $g \in G$. Since $\phi$ acts trivially on $G/N$, we have $g\phi = ng$ for some
$n \in N$. Since $\phi$ acts
trivially on $N$ it follows that for every $i$, $(g)\phi^{i} = n^{i}g$. In
particular, taking $i = |\phi|$ we find
that $n^{|\phi|}=1$. But $|\phi|$ is prime to $|G|$ so $n=1$, thus $g\phi = g$ as
desired.
\end{proof}

\begin{lem}
Let $G$ be a group, $p \in \pi(G)$ and $P$ a Sylow $p$-subgroup of $G$.
\\
(i) The number of conjugates of $P$ in $G$ is congruent to one modulo $p$.
\\
(ii) If $P \lhd G$ then $G$ splits over $P$.
\end{lem}


\begin{lem}
For any positive integer $n$ and prime power $q$, we have
\be
|{\hbox{GL}}(n,q)| = \prod_{i=0}^{n-1} (q^{n} - q^{i}).
\ee
\end{lem}


\begin{lem}
Let $G$ be a group of odd order. Then $G$ possesses an elementary abelian
characteristic $p$-subgroup, for
some $p \in \pi(G)$.
\end{lem}

\begin{proof}
A group is characteristically simple if and only if it is the internal direct
product of isomorphic simple
subgroups. Since every group of odd order is soluble, it follows that a
characteristically simple group
of odd order must be an elementary abelian $p$-group. This and the fact that
\be
B \; {\hbox{char}} \; A \; {\hbox{char}} \; G \; \Rightarrow \; B \; {\hbox{char}}
\; G,
\ee
allows us to prove the lemma by induction on the group order.
\end{proof}

\begin{lem}
If $\Omega(|{\hbox{Aut}}(G)|) \leq 4$ then $|{\hbox{Aut}}(G)|$ is even.
\end{lem}

\begin{proof}
This follows from the work of various authors. See the introduction to \cite{SZ}, where
all the necessary references are given.
\end{proof}

\begin{lem}
Let $p$ be an odd prime and $G$ a $p$-group such that Aut$(G)$ is also a $p$-group.
Then
$|{\hbox{Aut}}(G)| \geq 3^7$ if $p = 3$ and $|{\hbox{Aut}}(G)| \geq p^6$ if $p > 3$.
\end{lem}

\begin{proof}
See \cite{C1} and \cite{MS}.
\end{proof}

\begin{cor}
If $G$ is a group for which Aut$(G)$ is an odd order group of order strictly
less than $3^7$ then Aut$(G)$ is not a $p$-group.
\end{cor}

\begin{proof}
If Aut$(G)$ is a $p$-group then so is the central factor group $G/Z$, being
isomorphic to a normal subgroup of
the former. Thus $G$ is nilpotent, hence the internal direct product of its Sylow
subgroups, say $G = P_1
\times \cdots \times P_k$. But then also
\be
{\hbox{Aut}}(G) = \prod_{i=1}^{k} {\hbox{Aut}}(P_i).
\ee
Since the automorphism group of a $2$-group is either trivial or of even order, we
now obtain the desired
conclusion from Lemma 2.9.
\end{proof}

The next three results will be the main tools in our approach to proving Theorem
1.1. The first two already appear in several papers, so we omit proofs. See, for
example, Theorem 2.1 of \cite{C2}.

\begin{lem}
Suppose $G = AB$ for some pair of subgroups $A$ and $B$, where $A \lhd G$. 
Let $\phi \in {\hbox{Aut}}(A)$. Then $\phi$
extends to an automorphism $\phi^{*}$ of $G$ such that $|\phi| = |\phi^{*}|$
if the following two conditions are satisfied :
\par (i) $\phi$ is trivial on $A \cap B$,
\par (ii) in the group Aut$(A)$, $[\phi,A^{B}] = 1$.
\end{lem}


This lemma will be used in the case when $\phi$ has order 2 and the subgroup $A$ is
either abelian or a $p$-group of class 2. In the former case we have

\begin{cor}
Suppose $G = A \rtimes B$ where $A$ is abelian of order greater than 2. Then $G$
possesses an automorphism of order 2. In particular, if
$G$ has a normal, abelian Sylow $p$-subgroup, for some odd prime $p$, then 
$|{\hbox{Aut}}(G)|$ is even.
\end{cor}


In the case of non-abelian $A$ we have the following result. To the best of our
knowledge it is a new result, and we thus state it as a theorem in its own 
right. 

\begin{thm}
Let $G = AB$, where $A$ is a normal $p$-subgroup of class 2 with \\ $A/(Z \cap A)$
elementary abelian, $B/(Z \cap B)$ is an abelian $p^{\prime}$-group
and $A \cap B \subseteq Z$. Assume $B$ acts
non-trivially on $A$. Then $G$ possesses an
automorphism of order 2
.
\end{thm}

While our proof is lengthy, it is
crucial to our analysis in Section 3.

\begin{proof}
Since $B$ acts non-trivially on $A$, it does so even on $A/(Z \cap A)$, by Lemma 2.4.
Let $n$ be the rank of $A/(Z \cap A)$ as an elementary abelian $p$-group. The
group $A^{B}$ thus describes a non-trivial
representation of the abelian $p^{\prime}$-subgroup $B/(Z \cap B)$ in
GL$(n,p)$.
By Maschke's theorem, this representation is completely reducible, say
\be
A^B = \psi_1 \times \cdots \times \psi_r,
\ee
where rank$(\psi_i) = n_i$ say. Moreover,
Schur's lemma implies that $A^B/{\hbox{ker}} \psi_i$ is cyclic for each $i$.
We first prove
the theorem in two special cases, before establishing the general result.
\\
\\
{\sc Case 1} : Suppose that $A^B$ is irreducible.
\\
\\
Let $A^B = \; <I_b>$. Clearly, there is a choice of basis $Z\alpha_1,....,Z\alpha_n$
for $A/(Z \cap
A)$, an element $z \in Z$ and
non-negative integers $k_1,...,k_n$ such that the action of $b$ on $A$ is given by
\begin{eqnarray}
b^{-1} \alpha_i b = \alpha_{i+1}, \;\;\; i = 1,...,n-1, \\
b^{-1} \alpha_n b = z \alpha_{1}^{k_1} \cdots \alpha_{n}^{k_n}.
\end{eqnarray}
In other words, the matrix of $I_b$ with respect to this basis is
\be
M = M_b = \left( \begin{array}{ccccccc} 0 & 1 & 0 & \cdot & \cdot & \cdot & 0 \\ 0 &
0 & 1 & 0 & \cdot & \cdot & 0
\\ \cdot & \cdot & \cdot & \cdot & \cdot & \cdot & \cdot \\ \cdot & \cdot & \cdot &
\cdot & \cdot & \cdot &
\cdot \\ \cdot & \cdot & \cdot & \cdot & \cdot & \cdot & \cdot \\ 0 & 0 & \cdot &
\cdot & \cdot & 0 & 1 \\
k_1 & k_2 & \cdot & \cdot & \cdot & \cdot & k_n \end{array} \right).
\ee
Let $\mathbb{I}_n$ denote the $n \times n$ identity matrix. One checks easily that
det$(M - \mathbb{I}_n) = (-1)^{n+1}[(k_1 + k_2 + \cdots + k_n) - 1]$, hence, since
$M$ acts non-trivially and irreducibly, we have
\be
k_1 + k_2 + \cdots + k_n \not\equiv 1 \; ({\hbox{mod $p$}}).
\ee
Raising both sides of (2.7) to the $p$-th power for each $i$ yields
\be
\alpha_{1}^{p} = \alpha_{2}^{p} = \cdots = \alpha_{n}^{p}.
\ee
Doing the same for (2.8) and inserting (2.11) yields
\be
\alpha_{n}^{p} = z^{p} \alpha_{n}^{(k_1 + \cdots + k_n)p}.
\ee
By (2.10) there is a substitution $a_n := \sigma_n \alpha_n$, for some $\sigma_n \in
Z$, such that
$a_{n}^{p} = 1$. Working back through equations (2.7) one sees that one can 
make a sequence of similar
substitutions $a_i := \sigma_i \alpha_{i}$, $i = n-1,n-2,...,1$, such that
\be
a_{1}^{p} = a_{2}^{p} = \cdots = a_{n}^{p} = 1.
\ee
With respect to the elements $a_1,...,a_n$, there exist $z_1,...,z_n \in Z$ such
that the action of
$b$ is given by
\begin{eqnarray}
b^{-1}a_i b = z_i a_{i+1}, \;\;\; i = 1,...,n-1, \\
b^{-1} a_n b = z_n a_{1}^{k_1} \cdots a_{n}^{k_n}.
\end{eqnarray}
Now let $\zeta_1,...,\zeta_n$ be any elements of $Z$. By (2.13), the mapping
\be
z \mapsto z \; (z \in Z), \;\;\; a_i \mapsto \zeta_i a_{i}^{-1},
\ee
extends to a well-defined automorphism $\phi$ of $A$ of order $2$. Note that
$\phi$ is trivial on $Z \supseteq A \cap B$. Thus, by Lemma 2.11,
this extends to an
automorphism of $G$ if $[\phi,I_b] = 1$. From (2.14)-(2.16) this reduces to a system
of $n$
equations in the $n$ unknowns $\zeta_1,...,\zeta_n$, namely
\begin{eqnarray}
\zeta_i \zeta_{i+1}^{-1} = z_{i}^{2}, \;\;\; i = 1,...,n-1, \\
\zeta_{1}^{-k_1} \cdots \zeta_{n-1}^{-k_{n-1}} \zeta_{n}^{1-k_n} = z_{n}^{2}.
\end{eqnarray}
From (2.10) one deduces that this system has a unique solution. Thus $G$ possesses
an automorphism of
order $2$, as desired.
\\
\\
{\sc Case 2} : Suppose that 
$A^B = \psi_1 \times \psi_2$ where $\psi_2$ is one-dimensional
$(n_2 = 1)$ and trivial.
\\
\\
Let $\psi_1 = \; <I_b>$. From the previous case we can deduce that there is a basis
$Za_1,...,Za_{n-1},Za=Za_n$ for $A/(Z \cap A)$ such
that
\be
a_{1}^{p} = a_{2}^{p} = \cdots = a_{n-1}^{p} = 1
\ee
and that the action of $I_b$ is given by
\begin{eqnarray}
b^{-1}a_i b = z_i a_{i+1}, \;\;\; i = 1,...,n-2, \\
b^{-1} a_{n-1} b = z_{n-1} a_{1}^{k_1} \cdots a_{n-1}^{k_{n-1}}, \\
b^{-1}ab = a,
\end{eqnarray}
for some $z_1,...,z_{n-1},z \in Z$ and non-negative integers $k_1,...,k_{n-1}$
satisfying (2.10), with $n$ replaced by $n-1$.
The map $\phi : A \rightarrow
A$ should now satisfy
\be
z \mapsto z \; (z \in Z), \;\;\; a \mapsto a, \;\;\; a_i \mapsto \zeta_i a_{i}^{-1}
\;\; (i = 1,...,n-1).
\ee
We will have a well-defined automorphism of $A$ provided
\be
[a,a_i] = 1, \;\;\; {\hbox{for $i = 1,...,n-1$}}.
\ee
But (2.24) follows from (2.20)-(2.22) and (2.10). For inserting (2.22) into each of
the $n$ equations in (2.20) and (2.21) yields
\begin{eqnarray}
[a,a_i] = [a,a_{i+1}], \;\;\; i = 1,...,n-2, \\
\; [a,a_{n-1}] = \prod_{j=1}^{n-1} [a,a_j]^{k_j}.
\end{eqnarray}
Then (2.25) and (2.26), combined with (2.10), yields (2.24). The map $\phi$ extends
to $G$ for a unique choice of $\zeta_1,...,\zeta_{n-1}$, as in Case 1. This deals
with Case 2.
\\
\\
We now prove Theorem 2.13 in full generality. Suppose $A^B$ decomposes as in 
(2.6),
where $\psi_i = \; <I_{b_i}>$ acts irreducibly on the $n_i$-dimensional subspace $W_i$,
and suppose that $\psi_{1},...,\psi_{s}$ are non-trivial.
Then for each $i = 1,...,s$ there is a basis $Za_{i1},...,Za_{in_i}$ for $W_i$, and
bases $Za_1,...,Za_{r-s}$ for $W_{s+1},...,W_r$ respectively such that
\be
a_{ij}^{p} = 1, \;\;\; i = 1,...,s, \; j = 1,...,n_i,
\ee
and that the action of $B$ is fully described by
\begin{eqnarray}
b_{i}^{-1}a_{ij}b_i = z_{ij} a_{i,j+1}, \;\;\; i = 1,...,s, \; j = 1,...,n_i - 1, \\
b_{i}^{-1}a_{in_i}b = z_{in_i} a_{i1}^{k_{i1}} \cdots a_{in_i}^{k_{in_i}}, \;\;\; i =
1,...,s, \\
b_{i}^{-1}a_{lj} b_i = a_{lj}, \;\;\; {\hbox{whenever $i \neq l$}}, \\
b^{-1} a_m b = a_m, \;\;\; \forall \; b \in B, \; m = s+1,...,r,
\end{eqnarray}
for some $z_{ij} \in Z$ and non-negative integers $k_{ij}$ satisfying
\be
\sum_{j=1}^{n_i} k_{ij} \not\equiv 1 \; ({\hbox{mod $p$}}), \;\;\; i = 1,...,s.
\ee
Then the mapping
\be
z \mapsto z \; (z \in Z), \;\;\; a_{ij} \mapsto \zeta_{ij}
a_{ij}^{-1}, \;\;\; a_m \mapsto a_m,
\ee
extends, for a unique choice of the elements $\zeta_{ij} \in Z$, to an
automorphism of $A$ of order $2$, which acts trivially on $Z$ and
commutes with $A^B$. Thus $G$ has an automorphism of order $2$ and the
theorem is proved.
\end{proof}

\setcounter{equation}{0}

\section{Case-by-case analysis}

Throughout this section, $G$ denotes a potential counterexample to Theorem 1.1,
whose non-existence we shall establish. By Lemma 2.8, we have $\Omega(|{\hbox{Aut}}(G)|) \geq 5$ and, together with
Corollary 2.10, one readily checks that this leaves the following 19 possibilities
for $|{\hbox{Aut}}(G)|$ :
\begin{eqnarray}
3^5 \cdot p, \;\;\; p \in \{5,7\}, \\
3^4 \cdot p, \;\;\; p \in \{5,7,11,13,17,19,23\}, \\
5^4 \cdot 3, \\
3^4 \cdot 5^2, \\
3^3 \cdot p^2, \;\;\; p \in \{5,7\}, \\
5^3 \cdot 3^2, \\
3^3 \cdot 5 \cdot p, \;\;\; p \in \{7,11,13\}, \\
3^3 \cdot 7 \cdot 11, \\
3^2 \cdot 5^2 \cdot 7.
\end{eqnarray}
The order of $G/Z$ must divide one
of those on the above list. Let $p \in \pi(G/Z)$ and suppose that $p^2$ does not
divide $|G/Z|$. Then any
Sylow $p$-subgroup of $G$ is abelian. So if $G/Z$, and hence $G$, has a normal Sylow
$p$-subgroup, then Lemma
2.5(ii) and Corollary 2.12 imply that $G$ possesses an automorphism of order 2. Using
Lemmas 2.5(i) and 2.7 this fact, together with Lemma 2.9, is
easily checked to already rule out all but the following 21 possibilities for
$|G/Z|$ :
\begin{eqnarray}
3^i \cdot 5, \;\;\; i \in \{4,5\}, \\
3^i \cdot 13, \;\;\; i \in \{3,4\}, \\
3 \cdot 5^i, \;\;\; i \in \{2,3,4\}, \\
3^i \cdot 5^2, \;\;\; i \in \{2,3,4\}, \\
3^i \cdot 7^2, \;\;\; i \in \{1,2,3\}, \\
3^2 \cdot 5^3, \\
3^i \cdot 5 \cdot 7, \;\;\; i \in \{2,3\}, \\
3^i \cdot 5^2 \cdot 7, \;\;\; i \in \{1,2\}, \\
3^i \cdot 5 \cdot 11, \;\;\; i \in \{2,3\}, \\
3^3 \cdot 5 \cdot 13.
\end{eqnarray}
For these remaining possibilities we may thus assume that $G/Z$ contains no normal
subgroup of prime order. In what follows, there is no loss of generality in 
assuming that $\pi(G) = \pi(G/Z)$. For let $\pi$ be the set of primes dividing 
$|G|$ but not $|G/Z|$. By Lemma 2.5, we have $G = A \times G_1$, where $A$ is an
abelian Hall $\pi$-subgroup of $G$. Then, by Corollary 2.12, $G$ has an 
automorphism of order $2$ unless $A \cong C_1$ or $C_2$ and Aut$(G) =$ 
Aut$(G_1)$. \par
As our next step, we claim that in all but three of the above 21 cases, namely
\be
3^4 \cdot 5^2, \;\;\; 3^2 \cdot 5^2 \cdot 7, \;\;\; 3^4 \cdot 13,
\ee
we can identify a prime $p \in \pi(G/Z)$ such that the Sylow $p$-subgroup of $G/Z$
is normal of order $p^i$,
for some $i \in \{2,3,4,5\}$. More precisely we have the following table of values :
\\ \\

\begin{tabular}{|c|c|c|} \hline
$i$ & $p$ & $|G/Z|$ \\ \hline \hline
$2$ & $3$ & $3^2 \cdot 5 \cdot 7 \;\;$ $3^2 \cdot 5 \cdot 11$ \\ \hline
$2$ & $5$ & $3 \cdot 5^2 \;\;$ $3^2 \cdot 5^2 \;\;$ $3^3 \cdot 5^2 \;\;$ $3 \cdot
5^2 \cdot 7$ \\ \hline
$2$ & $7$ & $3 \cdot 7^2 \;\;$ $3^2 \cdot 7^2 \;\;$ $3^3 \cdot 7^2 \;\;$ \\ \hline
$3$ & $3$ & $3^3 \cdot 13 \;\;$ $3^3 \cdot 5 \cdot 7 \;\;$ $3^3 \cdot 5 \cdot 11\;\;$ $3^3 \cdot 5 \cdot 13$ \\
\hline
$3$ & $5$ & $3 \cdot 5^3 \;\;$ $3^2 \cdot 5^3$ \\ \hline
$4$ & $3$ & $3^4 \cdot 5$ \\ \hline
$4$ & $5$ & $3 \cdot 5^4$ \\ \hline
$5$ & $3$ & $3^5 \cdot 5$ \\ \hline
\end{tabular}

$\;$ \\
In 11 of the 18 cases listed in the table, the normality of the identified
$p$-subgroup is established by a direct
application of Lemma 2.5(i). In the case of $3^3 \cdot 13$ we also need to use our
assumption that $G_{13}$ is
not normal, being of prime order. The remaining 6 cases are those where $\omega(|G|)
= 3$. Here we appeal both to
the aforementioned assumption and to Lemma 2.7. One example will suffice to
illustrate the idea, say
$3^2 \cdot 5 \cdot 7$. By assumption neither $G_5$ nor $G_7$ is normal. Thus, by
Lemma 2.7, some $3$-subgroup
must be characteristic. If this is $G_3$, we are done, otherwise it is isomorphic to
$C_3$.
The quotient group then has order $3 \cdot 5 \cdot 7$. Applying Lemma 2.7 to the
quotient, we find that
$G/Z$ has a characteristic subgroup of order $9$, $15$ or $21$. In the first
instance we are done. In the
latter two, Lemma 2.5(i) and eq. (2.4) imply that $G_5$ resp. $G_7$ is characteristic
after all, contradicting our
assumptions.
\par Of the three numbers (3.20) omitted from the table, we can still apply Lemma 2.7
as above to conclude in the case of both $3^4 \cdot 5^2$ and $3^2 \cdot 5^2 \cdot 7$
that either the $G_3$ or $G_5$-subgroup must be normal,
though we don't know which one a priori. This leaves $3^4 \cdot 13$, where applying
the same kind of
analysis allows us only to conclude that either $G_3$ or a subgroup of
order $3^3$ must be normal. The case when $|G/Z| = 3^4 \cdot 13$ will be called the
{\em exceptional} case.
\\
\\
In the next part of our analysis, we ignore the exceptional case. This will be dealt
with at the finish.
In each of the remaining 20 cases, we are guaranteed to have a splitting 
$G = P \rtimes B$, where
\\
\\
(i) $P$ is a Sylow $p$-subgroup for some prime $p \in \{3,5,7\}$,
\\
(ii) $|P/(Z \cap P)| = p^i$ for some $i \in \{2,3,4,5\}$,
\\
(iii) $B$ is a $p^{\prime}$-group, either abelian or of class 2,
\\
(iv) If $p = 3$ then $\pi(B) \subseteq \{5,7,11,13\}$,
\\
(v) If $p = 5$ then $i \leq 4$ and $\pi(B) \subseteq \{3,7\}$,
\\
(vi) If $p = 7$ then $i=2$ and $\pi(B) = \{3\}$.
\\
\\
We now divide the analysis into 2 steps, according as to whether $P/(Z \cap P)$ is
elementary abelian or not.
\\
\\
{\sc Step 1} : $P/(Z \cap P)$ is elementary abelian.
\\
\\
If $B$ acts trivially on $P/(Z \cap P)$ then $G = P \times B$. It is easy to check,
using the tools of Section 2,
that in all cases $B$ possesses an automorphism of order $2$, hence so does $G$.
If $B$ acts non-trivially then all the conditions of Theorem 2.13 are satisfied, so
we are done here too.
\\
\\
{\sc Step 2} : $P/(Z \cap P)$ is not elementary abelian.
\\
\\
Then $P/(Z \cap P)$ must have order $p^3$ at least. In all but 4 cases, namely
\be
5^3 \cdot 3, \;\;\; 5^3 \cdot 3^2, \;\;\; 5^4 \cdot 3, \;\;\; 3^5 \cdot 5
\ee
we can apply Lemmas 2.4 and 2.6 directly to find that $B$ acts trivially on $P$, and 
then argue as in Step 1. We illustrate with
the example of $3^4 \cdot 5$. By Lemma 2.3(ii), there is a proper subgroup $P_2$ of
$P$, strictly containing
$Z \cap P$, and such that $P_{2}/(Z \cap P) \subseteq Z[P/(Z \cap P)]$. Thus $P_2$
is invariant under the action
of $B$, so $B$ acts on both $P_2/(Z \cap P)$ and
$P/P_2$, and both these groups are of order $3^i$ for some $i < 4$. Now Lemma 2.6
implies that
$|{\hbox{GL}}(i,3)|$ is not divisible by $5$ for any $i < 4$. It follows that $B$
acts trivially
on both $P_{2}/(Z \cap P)$ and $P/P_2$, and hence on $P/(Z \cap P)$ by Lemma 2.4.
\\
\\
This leaves us with the four cases in (3.21).
First suppose $|G/Z| \in \{5^3 \cdot 3, \; 5^3 \cdot 3^2\}$. In both cases, we
can choose a characteristic subgroup $P_2$ of $P$ as above such that $|P_2/(Z \cap
P)| = 5$. Then
$B$ centralises $P_2$, since $|{\hbox{GL}}(1,5)|$ is not divisible by $3$. But
$C_{G}(P_2)$
must then be a proper, characteristic subgroup of $G$ containing $P_2$ and $B$, and hence of
index
dividing $5^2$. Thus $G^{\prime} \subseteq C_{G}(P_2)$, which means that $B$ must
also act
trivially on $P/P_2$, thus on all of $P/(Z \cap P)$, by Lemma 2.4.
\par If $|G/Z| = 5^4 \cdot 3$ and $P/(Z \cap P)$ is non-abelian we can reach the
same conclusion,
for either the commutator subgroup or center of $P/(Z \cap P)$ must have order $5$
and yields
the characteristic subgroup playing the role of $P_2$ above. If $P/(Z \cap P)$ is
abelian, it must
be isomorphic to $C_{25} \times C_{25}$. In this case, we can still use the method of
Theorem 2.13. For one easily checks that $3 \parallel |{\hbox{Aut}}(C_{25} 
\times
C_{25})|$, and that if $B/(Z \cap B) := \; <Zb>$ acts non-trivially on $P/(Z \cap
P)$, then
there is a generating set $Za_1,Za_2$ for the latter such that
\be
b^{-1}a_1 b = a_2, \;\; b^{-1} a_2 b = z a_{1}^{-1} a_{2}^{-1}, \;\;\;\;\; \exists
\; z \in Z.
\ee
This was the starting point of the method of proving Case 1 of Theorem 2.13, and
the same argument
carries through here.
\par Finally, suppose $|P/(Z \cap P)| = 3^5 \cdot 5$. Let $P_{2}/(Z \cap P)$ be
a non-trivial
characteristic subgroup of $P/(Z \cap P)$ of index at least $3^2$, whose existence
is guaranteed
by Lemma 2.3(ii). As noted previously, $|{\hbox{GL}}(i,3)|$ is not divisible by $5$,
for any
$i < 4$. Hence we can deduce from Lemma 2.3 that $B$ acts trivially on $P$, unless
$|P_{2}/(Z \cap P)| = 5$ and $B$ acts irreducibly on $P/P_2$. But the latter is
impossible, since
$P_2 B$ is contained in the characteristic subgroup $C_{G}(P_2)$ of $G$. This completes
Step 2.
\\
\\
All that now remains in order to establish Theorem 1.1 is to handle the exceptional
case where $|G/Z| = 3^4 \cdot 13$. We can assume that the Sylow $3$-subgroup of $G$
is not
normal, as otherwise
Theorem 2.13 can be applied directly. As noted earlier, some $3$-subgroup $P$ such
that
$|P/(Z \cap P)| = 3^3$ must then be normal. Since $13$ does not divide
$|{\hbox{GL}}(2,3)|$,
then $P/(Z \cap P)$ must be elementary abelian and acted upon irreducibly by any
Sylow $13$-subgroup, or else the latter would also be normal in $G$. Moreover, the
quotient $G/P$ must be non-abelian of order $39$, isomorphic to $N_{G}(Q)/(Z \cap
Q)$, for
any Sylow $13$-subgroup $Q$ of $G$. These facts combined imply that
\par (i) $G/Z$ is a centerless group,
\par (ii) $|G/G^{\prime}|$ is a multiple of $3$,
\par (iii) $G = PB$ for some subgroup $B$ with $P \cap B \subseteq Z$.
\\
By Lemma 2.1, (i) implies that any non-trivial central automorphism of $G$ is outer.
But if $G$ possesses
any outer automorphisms, then $|{\hbox{Aut}}(G)|$ will have to be at
least $3^5 \cdot 13$, which is greater than $3^7$, a contradiction. Thus $G$ cannot
have any
non-trivial central
automorphisms. But then (ii) and Lemma 2.2 force the conclusion that $Z$ contains no
elements of order $3$. This and (iii) in turn imply that $G = P \rtimes B$ and that
$P$ is elementary abelian of order $27$. But now we can apply Corollary 2.12 to
conclude that $G$ has an automorphism of order 2. The theorem is proved.

\setcounter{equation}{0}

\section{Discussion}

A careful reading of the foregoing proof will show that we have established 
somewhat more than we have stated. In fact, if $G/Z$ has odd order less than 
$3^7$, then either $G$ has an automorphism
of order $2$, or $G/Z$ is a group of order $3^4 \cdot 13$ with a very specific
structure :
namely, it is the split extension of an elementary abelian group of order $3^3$ by
a non-abelian group of order $39$. In the latter case we could only show that either $G$ has
an automorphism of order $2$, or a non-trivial central automorphism which was
necessarily outer. It would be interesting to establish if this case is a 
genuine exception or whether a group with this structure cannot be N.I. either.
This might in turn shed light on the following questions :
\\
\\
1. What is
the smallest order of a non-nilpotent N.I. group and, conversely, 
the smallest order of a non-nilpotent automorphism group of odd order ? 
Theorem 1.1 and the result in \cite{MS} give the answers for nilpotent groups.
Martin's result \cite{M} that almost all $p$-groups have automorphism 
group a $p$-group also motivates giving special attention to non-nilpotent
N.I. groups.  
\\
\\
2. It should be particularly interesting to study complete groups, as these
have no non-trivial outer or central automorphisms. Several papers in 
the literature deal with the construction of complete groups of odd order
(see \cite{D}, \cite{H}, \cite{HR}, \cite{Sch} and \cite{So} for example), 
but the smallest possible order
of such a group is unknown. The smallest example in the
literature seems to be a complete group of order $3^{12} \cdot 5$ 
constructed in \cite{H}.
\\
\\
3. In \cite{SZ} the authors conjectured that if 
$\Omega (|{\hbox{Aut}}(G)|) \leq 
5$ then $|{\hbox{Aut}}(G)|$ is even (recall Lemma 2.8 above). This problem
remains open. Resolving it would also shed light on the previous two questions.
\\
\\ 
Finally, we remark that Corollary 2.12 and Theorem 2.13 are of some
interest in their own right. They give sufficient
conditions for an automorphism of order $2$ of a normal subgroup to be
extendible to the whole
group, when
the subgroup is nilpotent of class at most $2$. It would be interesting to 
find similar sufficient conditions when the subgroup
has class 3.

\section*{Acknowledgement}

We thank John Curran for helpful comments on an earlier draft and the referee 
for his suggestions. 
The research of the first author is partly supported by a grant from the
Swedish Research Council (Vetenskapsr\aa det).

\ \\

\end{document}